\documentclass{turing2012}
\pdfoutput=1

\usepackage{times}
\usepackage[pdftex]{graphicx}
\usepackage{latexsym}
\usepackage{tikz}
\usepackage{amsthm,amsfonts,amssymb,amsmath}

\def\defn{\emph}
\def\hat{\widehat}

\def\overM{{}\mskip3mu\overline{\mskip-3mu M\mskip-2mu}\mskip2mu}
\def\overR{{}\mskip3mu\overline{\mskip-3mu R\mskip-1mu}\mskip1mu}
\def\overS{{}\mskip3mu\overline{\mskip-3mu S\mskip-1mu}\mskip1mu}
\def\overU{{}\mskip2mu\overline{\mskip-2mu U\mskip-2mu}\mskip2mu}

\newtheorem{theorem}{Theorem}
\newtheorem{problem}{Problem}

\title{Evolution of unknotting strategies for knots and braids}
\author{Nicholas Jackson\institute{Mathematics Institute, University
    of Warwick, Coventry CV4 7AL, UK. Email: Nicholas.Jackson@warwick.ac.uk}
    \and Colin G. Johnson\institute{School of Computing, University of
    Kent, Canterbury CT2 7NF, UK. Email: C.G.Johnson@kent.ac.uk}}

\begin{document}

\maketitle

\begin{abstract}
This paper explores the problem of \emph{unknotting} closed braids and
classical knots in mathematical knot theory. We apply evolutionary
computation methods to learn sequences of moves that simplify knot
diagrams, and show that this can be effective both when the evolution
is carried out for individual knots and when a generic sequence of
moves is evolved for a set of knots.
\end{abstract}

\section{Introduction}

\subsection{Knots and links}
Knot theory is currently one of the richest and most vibrant areas of pure
mathematics, having connections not only with other topics in algebraic
and geometric topology, but also with many other branches of mathematics,
as well as mathematical physics~\cite{witten:qftjp} and
biochemistry~\cite{sumners:ltc}.

A full introduction to the study of knots and links is beyond the scope of
this article, but a readable introduction may be found in, for example,
the book by Cromwell~\cite{cromwell:ikl}, and a more comprehensive but
still accessible survey in Rolfsen's classic text~\cite{rolfsen:kl}.

\begin{figure}
\centering
\includegraphics[width=0.12\textwidth]{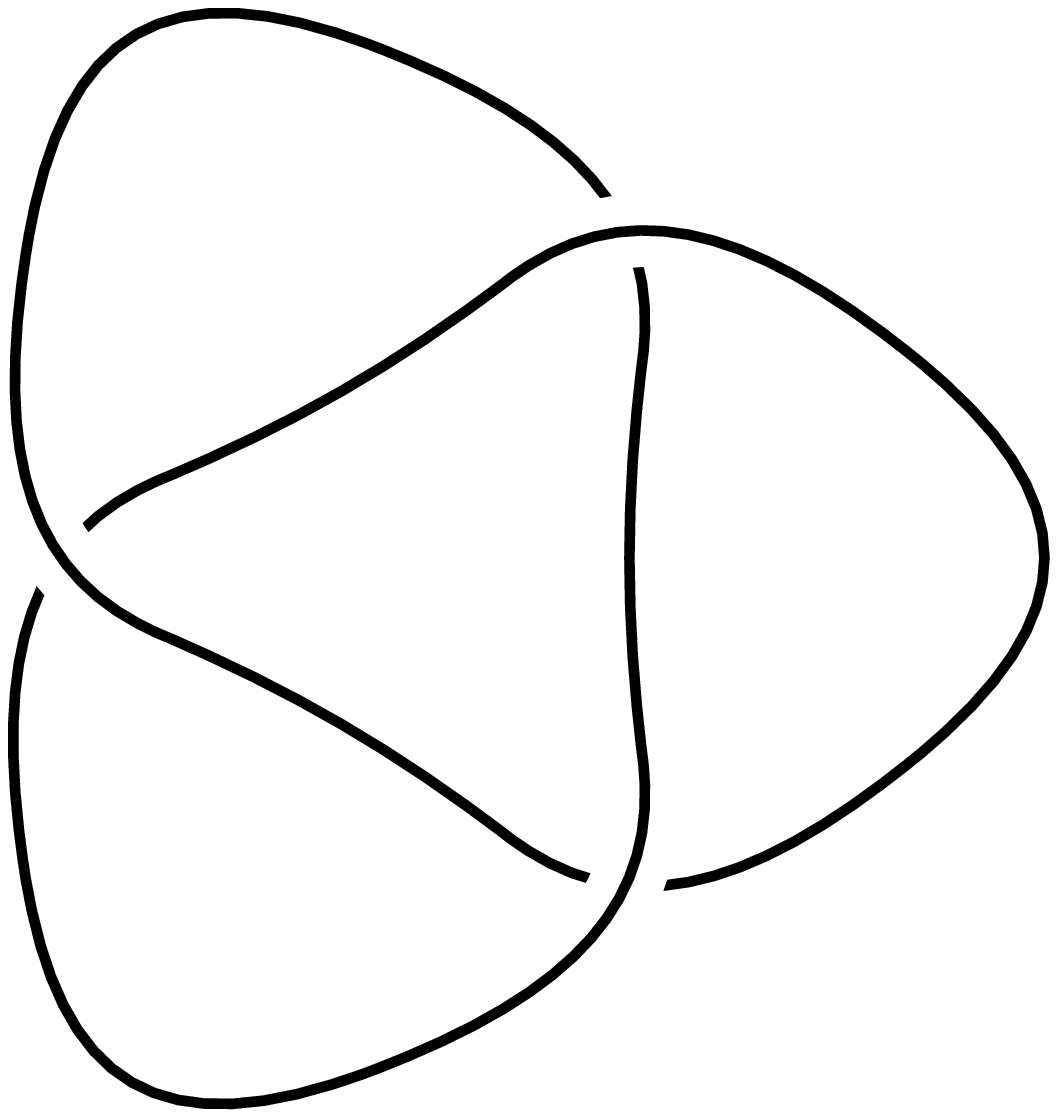}
\includegraphics[width=0.12\textwidth]{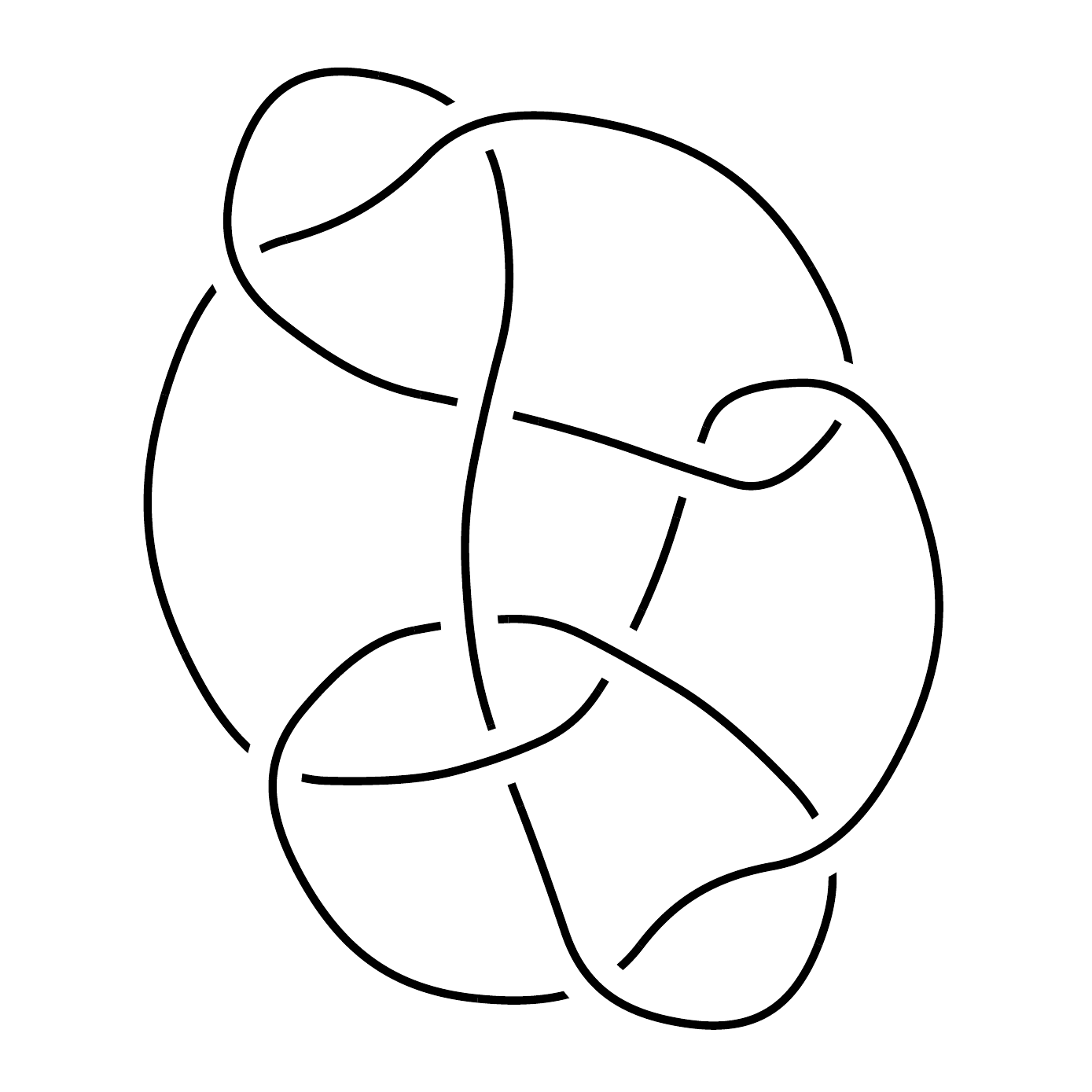}
\includegraphics[width=0.12\textwidth]{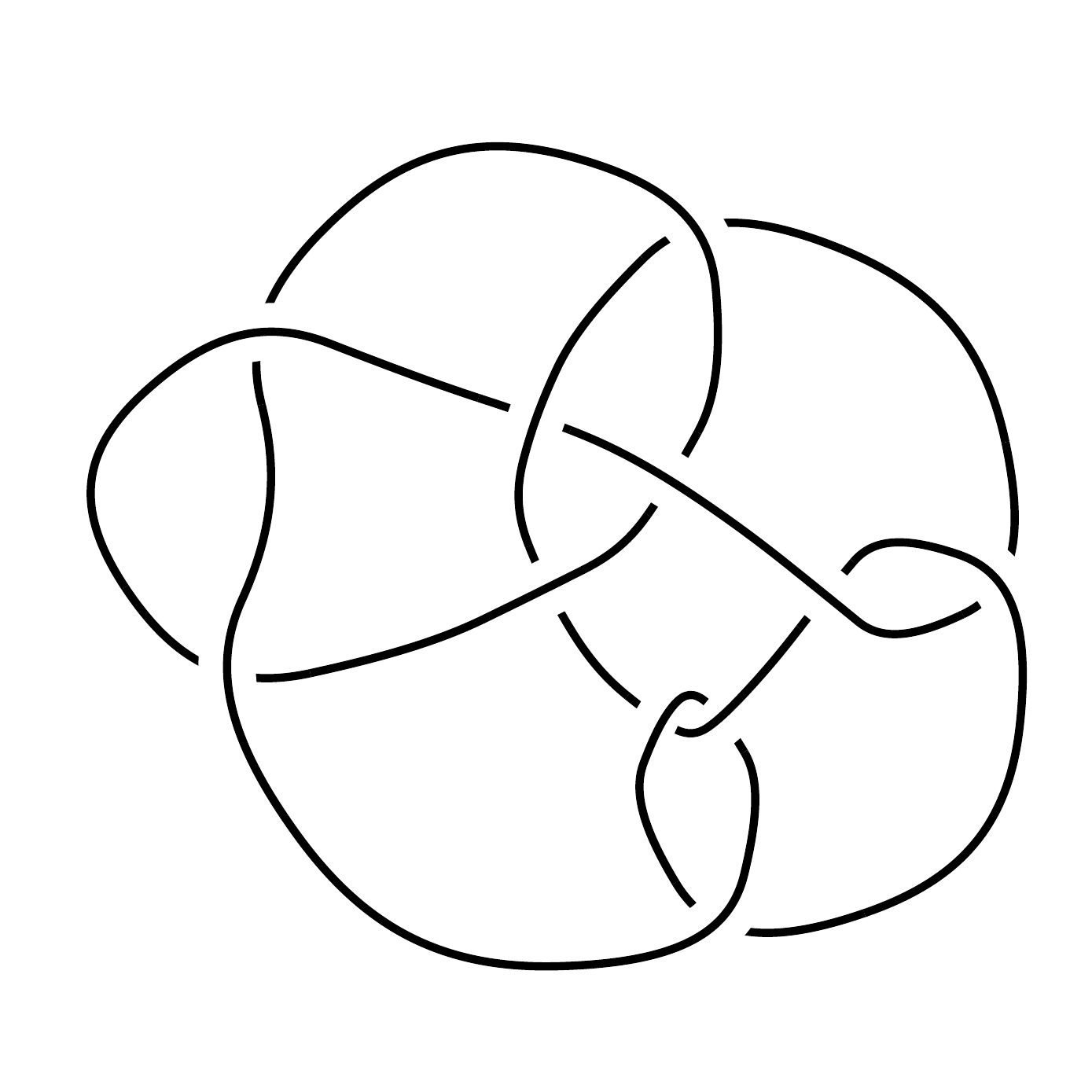}
\caption{Examples of knots: the trefoil ($3_1$), Conway's knot
($11n34$) and the Kinoshita--Terasaka knot ($11n42$)}
\label{fig:egknots}
\end{figure}

We define a \defn{knot} to be an isotopy class
of embeddings $K\colon S^1\hookrightarrow \mathbb{R}^3$, where $S^1 =
\{e^{i\theta} \in \mathbb{C} : 0\leqslant\theta<2\pi\}$ denotes the
standard unit circle; informally, this is a set of placements of a
closed circle in space.  A \defn{link} is a knot with more than one circular
component, that is, an (isotopy class of an) embedding
$L\colon S^1\sqcup\ldots\sqcup S^1\hookrightarrow\mathbb{R}^3$.
(Alternatively, a knot may be regarded as a link with a single component.)

We generally represent knots and links with diagrams in the plane:
projections of the embedded circle(s) where each double intersection point
is equipped with crossing information, and we disallow cusps, tangencies
or triple-points.  Examples may be seen in Figure~\ref{fig:egknots}.
Identifiers such as $3_1$ and $8_{17}$ refer to the table in Rolfsen's
book~\cite{rolfsen:kl}, while identifiers of the form $11a367$ and
$11n34$ refer to, respectively, alternating and non-alternating knots
in the census of Hoste, Thistlethwaite and Weeks~\cite{hoste:1701936}.

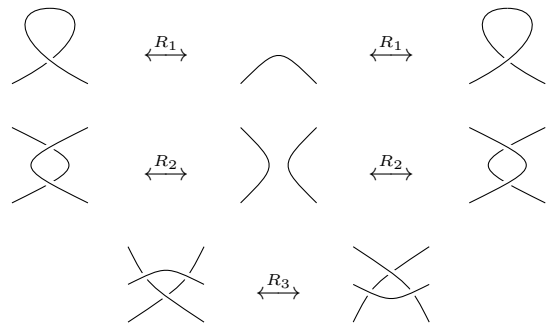
\begin{figure}
\centering
\begin{tabular}{ccccc}
\begin{tikzpicture}[scale=0.5,baseline={(0,0.3)}]
\draw (0,0) .. controls (2,1) and (2,2) .. (1,2);
\draw[white, line width=3pt] (1,2) .. controls (0,2) and (0,1) .. (2,0);
\draw (1,2) .. controls (0,2) and (0,1) .. (2,0);
\end{tikzpicture} &
$\stackrel{R_1}{\longleftrightarrow}$ &
\begin{tikzpicture}[scale=0.5,baseline={(0,0.3)}]
\draw (0,0) .. controls (1,1) .. (2,0);
\end{tikzpicture} &
$\stackrel{R_1}{\longleftrightarrow}$ &
\begin{tikzpicture}[scale=0.5,baseline={(0,0.3)}]
\draw (1,2) .. controls (0,2) and (0,1) .. (2,0);
\draw[white, line width=3pt] (0,0) .. controls (2,1) and (2,2) .. (1,2);
\draw (0,0) .. controls (2,1) and (2,2) .. (1,2);
\end{tikzpicture} \\[5ex]
\begin{tikzpicture}[scale=0.5,baseline={(0,0.3)}]
\draw (0,0) .. controls (2,1) .. (0,2);
\draw[white, line width=3pt] (2,0) .. controls (0,1) .. (2,2);
\draw (2,0) .. controls (0,1) .. (2,2);
\end{tikzpicture} &
$\stackrel{R_2}{\longleftrightarrow}$ &
\begin{tikzpicture}[scale=0.5,baseline={(0,0.3)}]
\draw (0,0) .. controls (1,1) .. (0,2);
\draw (2,0) .. controls (1,1) .. (2,2);
\end{tikzpicture} &
$\stackrel{R_2}{\longleftrightarrow}$ &
\begin{tikzpicture}[scale=0.5,baseline={(0,0.3)}]
\draw (2,0) .. controls (0,1) .. (2,2);
\draw[white, line width=3pt] (0,0) .. controls (2,1) .. (0,2);
\draw (0,0) .. controls (2,1) .. (0,2);
\end{tikzpicture} \\[5ex]
& \begin{tikzpicture}[scale=0.5,baseline={(0,0.3)}]
\draw (2,2) .. controls (1.5,1) .. (0,0);
\draw[white, line width=3pt] (0,2) .. controls (0.5,1) .. (2,0);
\draw (0,2) .. controls (0.5,1) .. (2,0);
\draw[white, line width=3pt] (0,1) .. controls (1,1.5) .. (2,1);
\draw (0,1) .. controls (1,1.5) .. (2,1);
\end{tikzpicture} &
$\stackrel{R_3}{\longleftrightarrow}$ &
\begin{tikzpicture}[scale=0.5,baseline={(0,0.3)}]
\draw (0,0) .. controls (0.5,1) .. (2,2);
\draw[white, line width=3pt] (2,0) .. controls (1.5,1) .. (0,2);
\draw (2,0) .. controls (1.5,1) .. (0,2);
\draw[white, line width=3pt] (0,1) .. controls (1,0.5) .. (2,1);
\draw (0,1) .. controls (1,0.5) .. (2,1);
\end{tikzpicture} &
\end{tabular}
\caption{Reidemeister moves}
\label{fig:rmoves}
\end{figure}

Isotopy of embeddings descends to certain allowable
local moves on diagrams which were first studied by
Reidemeister~\cite{reidemeister:ebk} and by
Alexander and Briggs~\cite{alexander:tkc}.  These \defn{Reidemeister
moves} are depicted in Figure~\ref{fig:rmoves}.  Two knots or links
are isotopic if and only if their diagrams are connected by a finite
sequence of Reidemeister moves and continuous deformations of the ambient
projection plane.

There are a number of different measures of the complexity of a given knot
or link $K$, the best known of which is the \defn{crossing number}: the
minimal number of crossings over all possible diagrams for $K$.  Related
to this is the \defn{unknotting number} $u(K)$: the minimal number, over
all possible diagrams for a knot, of crossings which must be changed in
order to obtain a trivial knot.  The trefoil in Figure~\ref{fig:egknots}
has unknotting number $u(3_1) = 1$: it may be seen to be nontrivially
knotted (that is, not isotopic to an unknotted circle) but changing any
single crossing results in a diagram which may be transformed (by means
of an $R_2$ move followed by an $R_1$ move) into an unknotted circle.
Less obviously, the other two knots in Figure~\ref{fig:egknots} are also
known to have unknotting number 1.

The unknotting number $u(K)$ is a conceptually simple measure of
the complexity of a given knot $K$ (broadly speaking, the higher the
unknotting number, the more knotted the knot in question) but one
which is often not straightforward to calculate.  According to the
\emph{KnotInfo} database~\cite{livingston:knotinfo}, the unknotting
number is currently unknown for nine of the 165 prime knots with
ten crossings (and also for many knots of higher crossing number); for
these nine knots, the unknotting number is is known to
be either 2 or 3, due at least in part to work on Heegaard Floer
homology by Oszv\'ath and Szab\'o~\cite{ozsvath-szabo:hfh-unknotting}
which ruled out unknotting number 1.  Recent work by Borodzik and
Friedl~\cite{borodzik-friedl:unknotting} has introduced a new invariant
which provides a lower bound of 3 for the unknotting number of twenty-five
otherwise difficult cases of knots with up to twelve crossings.

In this paper, our main goal is not necessarily to find optimal bounds on
the unknotting numbers of currently unresolved cases, but to explore the
possibilities afforded by applying evolutionary computing techniques to
pure mathematical problems in group theory and geometric topology, and to
try to obtain some qualitative understanding of the search landscape for
these problems.  The unknotting problem is relatively straightforward
to describe and implement, and thus provides a good candidate for a
preliminary investigation of this type.

\subsection{Braids}

The theory of braids was first seriously investigated by
Artin~\cite{artin:tdz}, and again a full treatment is far beyond the
scope of this article, so we will discuss only those aspects essential
for what follows.  A more comprehensive discussion may be found in
either the book by Hansen~\cite{hansen:bc} or the classic monograph by
Birman~\cite{birman:blmcg}.

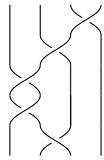
\begin{figure}
\centering
\begin{tikzpicture}[scale=0.4]
\draw (0,0) -- (0,1);
\draw (1,0) .. controls (1,0.4) and (2,0.6) .. (2,1);
\draw[white, line width=3pt] (2,0) .. controls (2,0.4) and (1,0.6) .. (1,1);
\draw (2,0) .. controls (2,0.4) and (1,0.6) .. (1,1);
\draw (3,0) -- (3,1);
\draw (1,1) .. controls (1,1.4) and (0,1.6) .. (0,2);
\draw[white, line width=3pt] (0,1) .. controls (0,1.4) and (1,1.6) .. (1,2);
\draw (0,1) .. controls (0,1.4) and (1,1.6) .. (1,2);
\draw (2,1) -- (2,2);
\draw (3,1) -- (3,2);
\draw (1,2) .. controls (1,2.4) and (0,2.6) .. (0,3);
\draw[white, line width=3pt] (0,2) .. controls (0,2.4) and (1,2.6) .. (1,3);
\draw (0,2) .. controls (0,2.4) and (1,2.6) .. (1,3);
\draw (2,1) -- (2,3);
\draw (3,2) -- (3,3);
\draw (0,3) -- (0,4);
\draw (2,3) .. controls (2,3.4) and (1,3.6) .. (1,4);
\draw[white, line width=3pt] (1,3) .. controls (1,3.4) and (2,3.6) .. (2,4);
\draw (1,3) .. controls (1,3.4) and (2,3.6) .. (2,4);
\draw (3,3) -- (3,4);
\draw (0,4) -- (0,5);
\draw (1,4) -- (1,5);
\draw (3,4) .. controls (3,4.4) and (2,4.6) .. (2,5);
\draw[white, line width=3pt] (2,4) .. controls (2,4.4) and (3,4.6) .. (3,5);
\draw (2,4) .. controls (2,4.4) and (3,4.6) .. (3,5);
\end{tikzpicture}
\caption{A 4--string braid}
\label{fig:braid}
\end{figure}

We define a geometric \defn{braid} on $n$ strings to be a system of $n$
disjoint, embedded arcs $A = \{A_1,\ldots,A_n\}$ in
$\mathbb{R}^2\times[0,1]$, such that the arc $A_i$ joins the point $P_i
= (i,0,1)$ to the point $Q_{\tau(i)} = (\tau(i),0,0)$, where $\tau$
denotes some permutation of the numbers $\{1,\ldots,n\}$, such that
each arc $A_i$ intersects the intermediate plane $\mathbb{R}^2\times\{z\}$
exactly once, for all $0<z<1$.  Figure~\ref{fig:braid} shows an example of
a 4--string braid.

The \defn{elementary braid} $\sigma_i$, for $1\leqslant i\leqslant n{-}1$,
is the $n$--string braid in which the $(i{+}1)$st string crosses over
the $i$th string, and no other interactions take place; its inverse
$\sigma_i^{-1}$ is the braid in which the $i$th string crosses over the
$(i{+}1)$st string (see Figure~\ref{fig:elbraid}).  Any $n$--string braid $\beta$ may be represented
(although not, in general, uniquely) as a concatenated sequence of
elementary braids.

\begin{figure}
\centering
\begin{tabular}{r}
$\sigma_i=$
\begin{tikzpicture}[scale=0.5,baseline={(0,0.2)}]
\draw (0,1) node [above] {\tiny 1} -- (0,0);
\draw (1,0.5) node {$\cdots$};
\draw (2,1) node [above] {\tiny $i{-}1$} -- (2,0);
\draw (3,1) node [above] {\tiny $i$} .. controls (3,0.6) and (4,0.4) .. (4,0);
\draw[white, line width=3pt] (4,1) .. controls (4,0.6) and (3,0.4) .. (3,0);
\draw (4,1) node [above] {\tiny $i{+}1$} .. controls (4,0.6) and (3,0.4) .. (3,0);
\draw (5,1) node [above] {\tiny $i{+}2$} -- (5,0);
\draw (6,0.5) node {$\cdots$};
\draw (7,1) node [above] {\tiny $n$} -- (7,0);
\end{tikzpicture}\\[3ex]
$\sigma_i^{-1}=$
\begin{tikzpicture}[scale=0.5,baseline={(0,0.2)}]
\draw (0,1) node [above] {\tiny 1}-- (0,0);
\draw (1,0.5) node {$\cdots$};
\draw (2,1) node [above] {\tiny $i{-}1$} -- (2,0);
\draw (4,1) node [above] {\tiny $i{+}1$} .. controls (4,0.6) and (3,0.4) .. (3,0);
\draw[white, line width=3pt] (3,1) .. controls (3,0.6) and (4,0.4) .. (4,0);
\draw (3,1) node [above] {\tiny $i$} .. controls (3,0.6) and (4,0.4) .. (4,0);
\draw (5,1) node [above] {\tiny $i{+}2$} -- (5,0);
\draw (6,0.5) node {$\cdots$};
\draw (7,1) node [above] {\tiny $n$} -- (7,0);
\end{tikzpicture}
\end{tabular}
\caption{The elementary braids $\sigma_i$ and $\sigma_i^{-1}$}
\label{fig:elbraid}
\end{figure}

We consider two $n$--braids $\beta_1$ and $\beta_2$ to be equivalent if
they are related by an isotopy which keeps their endpoints fixed.  In the
language of elementary braids, this translates to the following
identities:
\begin{align}
\sigma_i\sigma_j &= \sigma_j\sigma_i &&\text{for } |i-j|>1 \\
\sigma_i\sigma_{i+1}\sigma_i &= \sigma_{i+1}\sigma_i\sigma_{i+1}
&&\text{for } 1\leqslant i\leqslant n-1
\end{align}
Geometrically, the first of these corresponds to moving two
non-interacting elementary braids past each other, and the second is
essentially the third Reidemeister move $R_3$.

We may define the \defn{braid group} $B_n$ by the following presentation:
\begin{equation}
\bigg\langle \!\sigma_1,\ldots,\sigma_{n-1} \bigg| \!\!\begin{array}{ll}
  \sigma_i\sigma_j = \sigma_j\sigma_i & \!\!|i-j|>1 \\
  \sigma_i\sigma_{i+1}\sigma_i = \sigma_{i+1}\sigma_i\sigma_{i+1}
  & \!\!1\leqslant i\leqslant n-1 \end{array} \!\!\bigg\rangle
\end{equation}
It can be shown that the group defined by this presentation is isomorphic
to the group obtained by imposing the obvious concatenation operation on
the (in general, infinite) set of all $n$--string geometric braids.  In
this latter group, the identity element is the trivial $n$--braid (the one
with no crossings) and for any braid $\beta$, the inverse $\beta^{-1}$ may
be obtained by reflecting $\beta$ in the horizontal plane
$\mathbb{R}^2\times\big\{\frac12\big\}$.

There are other, equivalent constructions of the $n$--string braid
group, including one in terms of the fundamental group of a particular
configuration space factored by an action of the symmetric group $S_n$
but these will not concern us here.

\begin{figure}
\centering
\begin{tabular}{ccc}
\begin{tikzpicture}[scale=0.4,baseline={(0,0.7)}]
\draw (0,0) -- (0,1);
\draw (1,0) .. controls (1,0.4) and (2,0.6) .. (2,1);
\draw[white, line width=3pt] (2,0) .. controls (2,0.4) and (1,0.6) .. (1,1);
\draw (2,0) .. controls (2,0.4) and (1,0.6) .. (1,1);
\draw (1,1) .. controls (1,1.4) and (0,1.6) .. (0,2);
\draw[white, line width=3pt] (0,1) .. controls (0,1.4) and (1,1.6) .. (1,2);
\draw (0,1) .. controls (0,1.4) and (1,1.6) .. (1,2);
\draw (2,1) -- (2,2);
\draw (1,2) .. controls (1,2.4) and (0,2.6) .. (0,3);
\draw[white, line width=3pt] (0,2) .. controls (0,2.4) and (1,2.6) .. (1,3);
\draw (0,2) .. controls (0,2.4) and (1,2.6) .. (1,3);
\draw (2,1) -- (2,3);
\draw (0,3) -- (0,4);
\draw (1,3) .. controls (1,3.4) and (2,3.6) .. (2,4);
\draw[white, line width=3pt] (2,3) .. controls (2,3.4) and (1,3.6) .. (1,4);
\draw (2,3) .. controls (2,3.4) and (1,3.6) .. (1,4);
\end{tikzpicture}
& $\stackrel{\text{closure}}{\longrightarrow}$ &
\begin{tikzpicture}[scale=0.4,baseline={(0,0.7)}]
\draw (0,0) -- (0,1);
\draw (1,0) .. controls (1,0.4) and (2,0.6) .. (2,1);
\draw[white, line width=3pt] (2,0) .. controls (2,0.4) and (1,0.6) .. (1,1);
\draw (2,0) .. controls (2,0.4) and (1,0.6) .. (1,1);
\draw (1,1) .. controls (1,1.4) and (0,1.6) .. (0,2);
\draw[white, line width=3pt] (0,1) .. controls (0,1.4) and (1,1.6) .. (1,2);
\draw (0,1) .. controls (0,1.4) and (1,1.6) .. (1,2);
\draw (2,1) -- (2,2);
\draw (1,2) .. controls (1,2.4) and (0,2.6) .. (0,3);
\draw[white, line width=3pt] (0,2) .. controls (0,2.4) and (1,2.6) .. (1,3);
\draw (0,2) .. controls (0,2.4) and (1,2.6) .. (1,3);
\draw (2,1) -- (2,3);
\draw (0,3) -- (0,4);
\draw (1,3) .. controls (1,3.4) and (2,3.6) .. (2,4);
\draw[white, line width=3pt] (2,3) .. controls (2,3.4) and (1,3.6) .. (1,4);
\draw (2,3) .. controls (2,3.4) and (1,3.6) .. (1,4);
\draw (0,0) arc (180:360:3);
\draw (1,0) arc (180:360:2);
\draw (2,0) arc (180:360:1);
\draw (0,4) arc (180:0:3);
\draw (1,4) arc (180:0:2);
\draw (2,4) arc (180:0:1);
\draw (4,0) -- (4,4);
\draw (5,0) -- (5,4);
\draw (6,0) -- (6,4);
\end{tikzpicture} \\[15ex]
$\beta$ & $\stackrel{\text{closure}}{\longrightarrow}$ & $\hat\beta$
\end{tabular}
\caption{The closure operation on braids}
\label{fig:closure}
\end{figure}
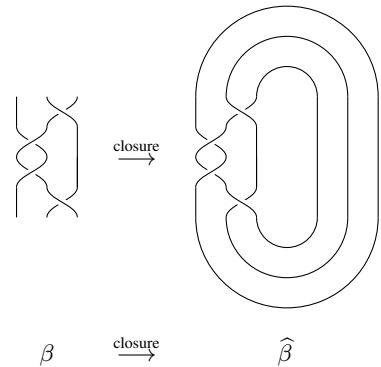

Given a braid $\beta\in B_n$, we can obtain a link $\hat\beta$ by
the \defn{closure} operation depicted in Figure~\ref{fig:closure},
that is, we join each point $P_i = (i,0,1)$ to the point $Q_i =
(i,0,0)$.  In fact, Alexander's Theorem~\cite{alexander:lskc} (see also
Birman~\cite[Theorem~2.1]{birman:blmcg}) states that any knot or link
can be obtained in this way; an explicit algorithm may be found in the
paper by Vogel~\cite{vogel:rlb}.  Note that a closed-braid presentation
need not be minimal with respect to the crossing number of the knot.  That
is, an $n$--crossing knot might not have a closed-braid presentation with
$n$ crossings.  Table~\ref{tbl:braidwords} lists several examples of
non-minimal presentations.

\begin{table}
\footnotesize
\centering
\begin{tabular}{|l|l|r|r|} \hline
$K$ & Braid word & Strands & Crossings \\ \hline
$3_1$ & $\sigma_{1}^3$ & 2 & 3 \\[1ex]
$4_1$ & $\sigma_{1} \sigma_{2}^{-1} \sigma_{1} \sigma_{2}^{-1}$ & 3 & 4
  \\[1ex]
$5_1$ & $\sigma_{1}^5$ & 2 & 5 \\
$5_2$ & $\sigma_{1}^3 \sigma_{2} \sigma_{1}^{-1} \sigma_{2}$ & 3 & 6
  \\[1ex]
$6_1$ & $\sigma_{1}^2 \sigma_{2} \sigma_{1}^{-1} \sigma_{3}^{-1}
  \sigma_{2} \sigma_{3}^{-1}$ & 4 & 7 \\
$6_2$ & $\sigma_{1}^3 \sigma_{2}^{-1} \sigma_{1} \sigma_{2}^{-1}$ & 3 & 6 \\
$6_3$ & $\sigma_{1}^2 \sigma_{2}^{-1} \sigma_{1} \sigma_{2}^{-2}$ & 3 & 6
  \\[1ex]
$7_1$ & $\sigma_{1}^7$ & 2 & 7 \\
$7_2$ & $\sigma_{1}^3 \sigma_{2} \sigma_{1}^{-1} \sigma_{2} \sigma_{3}
  \sigma_{2}^{-1} \sigma_{3}$ & 4 & 9 \\
$7_3$ & $\sigma_{1} \sigma_{1} \sigma_{1} \sigma_{1} \sigma_{1} \sigma_{2}
  \sigma_{1}^{-1} \sigma_{2}$ & 3 & 8 \\
$7_4$ & $\sigma_{1}^2 \sigma_{2} \sigma_{1}^{-1} \sigma_{2}^2 \sigma_{3}
  \sigma_{2}^{-1} \sigma_{3}$ & 4 & 9 \\
$7_5$ & $\sigma_{1}^4 \sigma_{2} \sigma_{1}^{-1} \sigma_{2}^2$ & 3 & 8 \\
$7_6$ & $\sigma_{1}^2 \sigma_{2}^{-1} \sigma_{1} \sigma_{3}
  \sigma_{2}^{-1} \sigma_{3}$ & 4 & 7 \\
$7_7$ & $\sigma_{1} \sigma_{2}^{-1} \sigma_{1} \sigma_{2}^{-1} \sigma_{3}
  \sigma_{2}^{-1} \sigma_{3}$ & 4 & 7 \\[1ex]
$8_1$ & $\sigma_{1}^2 \sigma_{2} \sigma_{1}^{-1} \sigma_{2}
  \sigma_{3} \sigma_{2}^{-1} \sigma_{4}^{-1} \sigma_{3} \sigma_{4}^{-1}$
  & 5 & 10 \\
$8_2$ & $\sigma_{1}^5 \sigma_{2}^{-1} \sigma_{1} \sigma_{2}^{-1}$ & 3 & 8 \\
$8_3$ & $\sigma_{1}^2 \sigma_{2} \sigma_{1}^{-1} \sigma_{3}^{-1}
  \sigma_{2} \sigma_{3}^{-1} \sigma_{4}^{-1} \sigma_{3} \sigma_{4}^{-1}$
  & 5 & 10 \\
$8_4$ & $\sigma_{1}^3 \sigma_{2}^{-1} \sigma_{1}
  \sigma_{2}^{-1} \sigma_{3}^{-1} \sigma_{2} \sigma_{3}^{-1}$ & 4 & 9 \\
$8_5$ & $\sigma_{1}^3 \sigma_{2}^{-1} \sigma_{1}^3 \sigma_{2}^{-1}$
  & 3 & 8 \\
$8_6$ & $\sigma_{1}^4 \sigma_{2} \sigma_{1}^{-1} \sigma_{3}^{-1}
  \sigma_{2} \sigma_{3}^{-1}$ & 4 & 9 \\
$8_7$ & $\sigma_{1}^4 \sigma_{2}^{-1} \sigma_{1} \sigma_{2}^{-1}
  \sigma_{2}^{-1}$ & 3 & 8 \\
$8_8$ & $\sigma_{1}^3 \sigma_{2} \sigma_{1}^{-1} \sigma_{3}^{-1}
  \sigma_{2} \sigma_{3}^{-2}$ & 4 & 9 \\
$8_9$ & $\sigma_{1}^3 \sigma_{2}^{-1} \sigma_{1} \sigma_{2}^{-3}$ & 3 & 8 \\
$8_{10}$ & $\sigma_{1}^3 \sigma_{2}^{-1} \sigma_{1}^2 \sigma_{2}^{-2}$
  & 3 & 8 \\
$8_{11}$ & $\sigma_{1}^2 \sigma_{2} \sigma_{1}^{-1} \sigma_{2}^2
  \sigma_{3}^{-1} \sigma_{2} \sigma_{3}^{-1}$ & 4 & 9 \\
$8_{12}$ & $\sigma_{1} \sigma_{2}^{-1} \sigma_{1} \sigma_{3} \sigma_{2}^{-1}
  \sigma_{4}^{-1} \sigma_{3} \sigma_{4}^{-1}$ & 5 & 8 \\
$8_{13}$ & $\sigma_{1}^2 \sigma_{2}^{-1} \sigma_{1} \sigma_{2}^{-2}
  \sigma_{3}^{-1} \sigma_{2} \sigma_{3}^{-1}$ & 4 & 9 \\
$8_{14}$ & $\sigma_{1}^3 \sigma_{2} \sigma_{1}^{-1}
  \sigma_{2} \sigma_{3}^{-1} \sigma_{2} \sigma_{3}^{-1}$ & 4 & 9 \\
$8_{15}$ & $\sigma_{1}^2 \sigma_{2}^{-1} \sigma_{1} \sigma_{3}
  \sigma_{2}^3 \sigma_{3}$ & 4 & 9 \\
$8_{16}$ & $\sigma_{1}^2 \sigma_{2}^{-1} \sigma_{1}^2 \sigma_{2}^{-1}
  \sigma_{1} \sigma_{2}^{-1}$ & 3 & 8 \\
$8_{17}$ & $\sigma_{1}^2 \sigma_{2}^{-1} \sigma_{1} \sigma_{2}^{-1}
  \sigma_{1} \sigma_{2}^{-2}$ & 3 & 8 \\
$8_{18}$ & $\sigma_{1} \sigma_{2}^{-1} \sigma_{1} \sigma_{2}^{-1} \sigma_{1}
  \sigma_{2}^{-1} \sigma_{1} \sigma_{2}^{-1}$ & 3 & 8 \\
$8_{19}$ & $\sigma_{1}^3 \sigma_{2} \sigma_{1}^3 \sigma_{2}$ & 3 & 8 \\
$8_{20}$ & $\sigma_{1}^3 \sigma_{2}^{-1} \sigma_{1}^{-3} \sigma_{2}^{-1}$
  & 3 & 8 \\
$8_{21}$ & $\sigma_{1}^3 \sigma_{2} \sigma_{1}^{-2} \sigma_{2}^2$
  & 3 & 8 \\ \hline
\end{tabular}
\caption{Braid words for knots with up to eight
crossings~\cite{livingston:knotinfo}}
\label{tbl:braidwords}
\end{table}

The following theorem gives explicit conditions
for when two different braids yield isotopic knots or links.
This result was due originally to Markov, although
a full proof was only published some years later by
Birman~\cite[Theorem~2.3]{birman:blmcg} (see also the paper by Rourke
and Lambropoulou~\cite{lambropoulou:mt3m}).

\begin{figure}
\centering
\begin{tabular}{ccccc}
& \begin{tikzpicture}[scale=0.35,baseline={(0,1)}]
\draw (0,0) -- (3,0) -- (3,2) -- (0,2) -- (0,0);
\draw (0,3) -- (3,3) -- (3,5) -- (0,5) -- (0,3);
\draw (0.5,2) -- (0.5,3); \draw (2.5,2) -- (2.5,3);
\draw (1.6,2.5) node {$\cdots$};
\draw (0.5,0) -- (0.5,-1); \draw (2.5,0) -- (2.5,-1);
\draw (1.6,-0.5) node {$\cdots$};
\draw (0.5,5) -- (0.5,6); \draw (2.5,5) -- (2.5,6);
\draw (1.6,5.5) node {$\cdots$};
\draw (1.6,1) node {$\beta_2$};
\draw (1.6,4) node {$\beta_1$};
\end{tikzpicture}
& $\stackrel{M_1}{\longleftrightarrow}$ &
\begin{tikzpicture}[scale=0.35,baseline={(0,1)}]
\draw (0,0) -- (3,0) -- (3,2) -- (0,2) -- (0,0);
\draw (0,3) -- (3,3) -- (3,5) -- (0,5) -- (0,3);
\draw (0.5,2) -- (0.5,3); \draw (2.5,2) -- (2.5,3);
\draw (1.6,2.5) node {$\cdots$};
\draw (0.5,0) -- (0.5,-1); \draw (2.5,0) -- (2.5,-1);
\draw (1.6,-0.5) node {$\cdots$};
\draw (0.5,5) -- (0.5,6); \draw (2.5,5) -- (2.5,6);
\draw (1.6,5.5) node {$\cdots$};
\draw (1.6,1) node {$\beta_1$};
\draw (1.6,4) node {$\beta_2$};
\end{tikzpicture} & \\[11ex]
\begin{tikzpicture}[scale=0.35,baseline={(0,0.7)}]
\draw (0,1) -- (3,1) -- (3,3) -- (0,3) -- (0,1);
\draw (0.5,1) -- (0.5,0); \draw (0.5,3) -- (0.5,4); \draw (2.5,3) -- (2.5,4);
\draw (3.5,1) -- (3.5,4);
\draw (2.5,1) .. controls (2.5,0.6) and (3.5,0.4) .. (3.5,0);
\draw[white, line width=3pt] (3.5,1) .. controls (3.5,0.6) and (2.5,0.4) .. (2.5,0);
\draw (3.5,1) .. controls (3.5,0.6) and (2.5,0.4) .. (2.5,0);
\draw (1.6,2) node {$\beta$};
\draw (0.5,4) node[above] {\tiny 1};
\draw (2.5,4) node[above] {\tiny $n$};
\draw (3.5,4) node[above] {\tiny $n{+}1$};
\draw (1.6,3.5) node {$\cdots$};
\draw (1.6,0.5) node {$\cdots$};
\end{tikzpicture}
& $\stackrel{M_2}{\longleftrightarrow}$ &
\begin{tikzpicture}[scale=0.35,baseline={(0,0.7)}]
\draw (0,1) -- (3,1) -- (3,3) -- (0,3) -- (0,1);
\draw (0.5,1) -- (0.5,0); \draw (0.5,3) -- (0.5,4);
\draw (2.5,1) -- (2.5,0); \draw (2.5,3) -- (2.5,4);
\draw (1.6,2) node {$\beta$};
\draw (0.5,4) node[above] {\tiny 1};
\draw (2.5,4) node[above] {\tiny $n$};
\draw (1.6,3.5) node {$\cdots$};
\draw (1.6,0.5) node {$\cdots$};
\end{tikzpicture}
& $\stackrel{M_2}{\longleftrightarrow}$ &
\begin{tikzpicture}[scale=0.35,baseline={(0,0.7)}]
\draw (0,1) -- (3,1) -- (3,3) -- (0,3) -- (0,1);
\draw (0.5,1) -- (0.5,0); \draw (0.5,3) -- (0.5,4); \draw (2.5,3) -- (2.5,4);
\draw (3.5,1) -- (3.5,4);
\draw (3.5,1) .. controls (3.5,0.6) and (2.5,0.4) .. (2.5,0);
\draw[white, line width=3pt] (2.5,1) .. controls (2.5,0.6) and (3.5,0.4) .. (3.5,0);
\draw (2.5,1) .. controls (2.5,0.6) and (3.5,0.4) .. (3.5,0);
\draw (1.6,2) node {$\beta$};
\draw (0.5,4) node[above] {\tiny 1};
\draw (2.5,4) node[above] {\tiny $n$};
\draw (3.5,4) node[above] {\tiny $n{+}1$};
\draw (1.6,3.5) node {$\cdots$};
\draw (1.6,0.5) node {$\cdots$};
\end{tikzpicture}
\end{tabular}
\caption{Markov moves of type $M_1$ (conjugation) and $M_2$
(stabilisation)}
\label{fig:markov}
\end{figure}
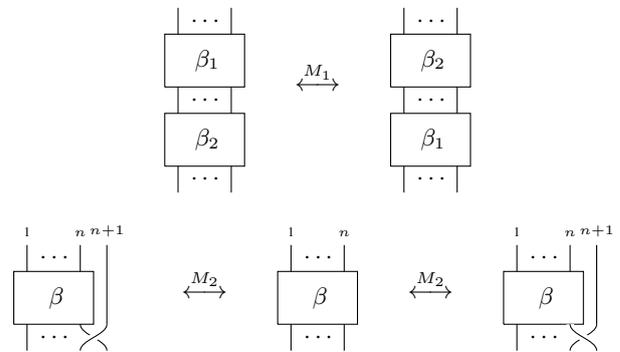

\begin{theorem}[Markov~\cite{markov:fagz}]
Two braids $\beta_1\in B_m$ and $\beta_2\in B_n$ yield closures
$\hat\beta_1$ and $\hat\beta_2$ which are isotopic as links if and only if
$\beta_1$ and $\beta_2$ are connected by a finite sequence of moves of
type $M_1$ (conjugation) and $M_2$ (stabilisation), as depicted in
Figure~\ref{fig:markov}.
\end{theorem}

\section{Unknotting}

We now attempt to use evolutionary techniques to devise optimal unknotting
strategies for knots and links.  The theorems of Alexander and Markov
enable us to represent knots and links as words in the standard generators
$\sigma_i$ of the braid group $B_n$ for some $n$.  The crossing-change
operation is then simply a matter of taking such a word and then replacing
a given $\sigma_i$ with its inverse $\sigma_i^{-1}$ or vice-versa.

Our goal is, given a knot $K$ represented as the closure of a braid word
$w\in B_n$, to evolve a sequence of certain moves which trivialises
the knot with the smallest number of crossing-changes, thus obtaining
an upper bound on the unknotting number $u(K)$ of $K$.  The allowable
moves are those which either leave the isotopy class of the corresponding
knot unchanged, or which change the sign of a single crossing, and may
be seen in Table~\ref{tbl:moves}.

\begin{table}
\centering
\begin{tabular}{|l|rcl|} \hline
$R_2^{\pm}$ & $\sigma_i^{\pm1}\sigma_i^{\mp1}$ & $\longmapsto$ & $1$\\
$\overR_2^{\pm}$ & $1$ & $\longmapsto$ &
$\sigma_i^{\pm1}\sigma_i^{\mp1}$ \\[1ex]
$R_3^{\pm}$ & $\sigma_i^{\pm1}\sigma_{i+1}^{\pm1}\sigma_{i}^{\pm1}$ & $\longmapsto$ &
  $\sigma_{i+1}^{\pm1}\sigma_i^{\pm1}\sigma_{i+1}^{\pm1}$ \\
$\overR_3^{\pm}$ & $\sigma_{i+1}^{\pm1}\sigma_i^{\pm1}\sigma_{i+1}^{\pm1}$ &
$\longmapsto$ & $\sigma_i^{\pm1}\sigma_{i+1}^{\pm1}\sigma_{i}^{\pm1}$ \\[1ex]
$M_1^{\pm}$ & $\sigma_i^{\pm1}\alpha$ & $\longmapsto$ &
  $\alpha\sigma_i^{\pm1}$ \\
$\overM_1^{\pm}$ & $\alpha\sigma_i^{\pm1}$ & $\longmapsto$ & $\sigma_i^{\pm1}\alpha$ \\[1ex]
$M_2^{\pm}$ & $\alpha\sigma_n^{\pm1}$ & $\longmapsto$ & $\alpha$ \\ 
$\overM_2^{\pm}$ & $\alpha$ & $\longmapsto$ & $\alpha\sigma_n^{\pm1}$
\\[1ex]
$S^{\pm\pm}$ &
  $\sigma_i^{\pm1}\sigma_j^{\pm1}$ & $\longmapsto$ &
  $\sigma_j^{\pm1}\sigma_i^{\pm1}$ (if $|i{-}j|>1$)\\
$\overS^{\pm\pm}$ &
  $\sigma_j^{\pm1}\sigma_i^{\pm1}$ & $\longmapsto$ &
  $\sigma_i^{\pm1}\sigma_j^{\pm1}$ (if $|i{-}j|>1$)\\[1ex]
$U^+ = \overU^-$ & $\sigma_i$ & $\longmapsto$ & $\sigma_i^{-1}$ \\
$U^- = \overU^+$ & $\sigma_i^{-1}$ & $\longmapsto$ & $\sigma_i$ \\
\hline
\end{tabular}
\caption{Allowable moves on braid words}
\label{tbl:moves}
\end{table}

In more detail, a move of type $R_2^+$ cancels a substring of the form
$\sigma_i\sigma_i^{-1}$, a move of type $R_2^-$ cancels a substring
of the form $\sigma_i^{-1}\sigma_i$, and moves of type $\overR^{\pm}$
introduce corresponding substrings.  These are the Reidemeister moves of
type $2$, translated into the context of braid words.  Similarly, moves
of type $R_3$ perform Reidemeister moves of type $3$, moves of type $M_1$
and $M_2$ are Markov moves, moves of type $S$ represent non-interacting
crossings sliding past each other, and moves of type $U$ change the sign
of a single crossing.

\section{Methods}\label{se:methods}

In this paper, we use evolutionary techniques to find sequences of
unknotting primitives which are optimal with respect to two subtly
different but related problems:

\begin{problem}\label{prob:single}
Given an arbitrary knot $K$, described as the closure $\hat{\beta}$
of some braid word $\beta\in B_n$, is there a sequence of moves which
reduces $\beta$ to the trivial word $1\in B_1$?  If so, what is the
minimal such sequence with respect to the number of crossing-change
operations $U^{\pm}$, and in what cases does this yield a sharp upper
bound for the unknotting number $u(K)$ of $K$ when compared
to known values such as those listed in the  \emph{KnotInfo}
database~\emph{\cite{livingston:knotinfo}}.
\end{problem}

\begin{problem}\label{prob:multiple}
Given a finite set $S$ of knots, each described as the closure
$\hat{\beta}$ of some braid word $\beta\in B_n$, is there a single,
universal sequence of moves which, possibly with repeated applications,
trivialises each knot in $S$?
\end{problem}

The second of these problems is, broadly speaking, the generalisation of
the first to more than one reference knot; equivalently, the first problem
is the special case of the second where we consider a single knot.

As usual in an evolutionary computation approach, we need to define
seven things: how population members are represented; the parameters;
how the population is initialised; how the mutation and crossover
operators are defined; and how fitness evaluation and selection
happen.

\subsection{Representation}

Each member of the population consists of a list of allowed moves drawn
from those in Table~\ref{tbl:moves}. This list is of variable length.

\subsection{Parameters}

For Problem~\ref{prob:single} the population size was 500 and the
number of iterations was $4\times \mathrm{length}(K)^2$, where $K$ is
the knot in question and length is the number of crossings. For
Problem~\ref{prob:multiple} the population size was 200 and the number
of iterations was $4\times \max(\mathrm{length}(S))^2$, where $S$ is
the set of knots. These parameters were determined by informal
experimentation. The program was run three times for each knot or set
of knots, and the best result is reported.

\subsection{Initialisation}

Each member of the population is initialised by selecting a list
of between 1 and 15 moves uniformly at random from the moves given
in Table~\ref{tbl:moves}.  The moves are selected with replacement;
that is, a sequence may include more than one move of a given type.

\subsection{Mutation}
There are three different mutation
operators, which are selected uniformly at random and applied to
individuals with the overall mutation probability being 10\%. These are:
\begin{enumerate}
\item Select a random move from the list and replace it with a move
drawn randomly from the moves in Table~\ref{tbl:moves}.
\item Choose a random move from the list and delete it (as long as the
list of moves contains at least one move).
\item Choose a random position in the list and insert a
randomly chosen move from the moves in Table~\ref{tbl:moves}.
\end{enumerate}

\subsection{Crossover}

One-point crossover is applied to all individuals as follows: the two
strings are aligned, a position less than or equal to the length of the
shortest string is selected at random, and the strings crossed over at
that point.

\subsection{Fitness evaluation}

The members of the population are evaluated by attempting to unknot each
of a set of knots, which in the case of Problem~\ref{prob:single} will
be just a single example, and in the case of Problem~\ref{prob:multiple}
will consist of more than one knot.  The execution is carried out as
follows. Let $M={m_1,m_2,\ldots,m_n}$ be the $n$ moves in the list. Let
$K = K_0$ be the original knot, and $K_1, K_2,\ldots$ be the sequence
of knots generated.

The knot $K_0$ is analysed for the preconditions for $m_1$ to be
carried out, if they are satisfied then the move is applied, so that
$K_1:=m_1(K_0)$; if the preconditions are not satisfied then the knot
is unchanged ($K_1:=K_0$). The next step is to attempt to apply $m_2$
to $K_1$ by seeing if its preconditions are satisfied, and so on. When
the end of $M$ is reached, the list is begun again from the beginning.

Each time the knot is changed as a result of applying the move, the
knot is checked to see if there are any crossings remaining. If so,
the algorithm terminates, and a positive result returned.

For Problem~\ref{prob:single} we apply the sequence $M$ once only,
but for Problem~\ref{prob:multiple} we perform repeated applications
of $M$. If the knot hasn't been trivialised by 50 applications of $M$,
we assume that it has become stuck in a repeating loop (which, for the
1\,701\,936 knots with sixteen or fewer crossings~\cite{hoste:1701936},
is a valid assumption), terminate the process with a negative result,
and move on to the next knot (if any) on the list.

By inspection of Table~\ref{tbl:moves} we see that certain operations
(namely, those of type $R_2$ and $M_2$) reduce the length of the braid word
under investigation, some (types $R_3$, $M_1$ and $U$) don't, and the rest
(types $\overR_2$ and $\overM_2$) increase the length of the braid word.

It is known (see, for example, the paper by Coward~\cite{coward:ordering})
that there exist diagrams for the unknot which can only be reduced
to the standard (zero-crossing) diagram of the unknot by means of at
least one move of type $\overR_2$.  That is, at some point during the
reduction process, we have to temporarily increase the complexity.
In this article however, for simplicity, we will restrict ourselves to
sequences of moves which do not increase the length of the braid word.
We intend to explore the more general case in later work.

With regard to Problem~\ref{prob:single}, we are trying to find a
sequence of operations which unknot a specific knot with the smallest
number of crossing changes.  Slightly less importantly, we want to find
the simplest possible such unknotting sequence.

For Problem~\ref{prob:multiple}, we are trying to find a sequence of
operations which (perhaps with repeated applications) unknots as many
knots as possible, as efficiently as possible.

Let $S$ denote the set of reference braids, on which each sequence is
being tested.  (For Problem~\ref{prob:single} this will consist of a
single braid word.)  Let $r_S(M)$ denote the number of braids in $S$
which are fully reduced by (one or more application of) the sequence
$M$.  (In practice, we specify an upper threshold of 50 iterations, as
described above.)  Let $\min_S(M)$ and $\max_S(M)$ denote, respectively,
the minimum and maximum number, over all braids in $S$, of iterations
of $M$ required to reduce a (reducible) braid.  Let $l(M)$ denote
the length of the sequence $M$.  By $c(M)$ we denote the number of
crossing-change (type $U$) operations in $M$, and by $c_S(M)$ we denote
the total number, over all braids in $S$, of successful crossing-change
(type $U$) operations.  That is, $c_S(M)$ gives a measure of the total
amount of unknotting actually performed by the sequence $M$.

In the case of Problem~\ref{prob:single}, some of the operations in
the string $M$ may have no effect on the braid under examination.
For example, applying an $R_2$ move to a braid which at that stage has
no $\sigma_i^{\pm1}\sigma_i^{\mp1}$ substrings will leave the braid
unchanged, and may thus be safely elided from the sequence, resulting in
a shorter sequence.  Given a sequence $M$, applied to a braid $\beta\in
B_n$, we denote by $l_{\text{opt}}(M)$ the length of the sequence obtained
by optimising $M$ in this way with respect to $\beta$.

The fitness function should, ideally, seek to minimise the
number of crossing changes, maximise (at least when working on
Problem~\ref{prob:multiple}) the number of knots which can be reduced by
a given sequence, minimise the length of the sequence, and minimise the
number of repeated applications of the sequence required to reduce those
braids in $S$ which are reducible by the operations under consideration.

With those criteria in mind, we define the fitness function for
Problem~\ref{prob:single} to be
$$f_1(M) = 1 + \frac{10000 r_S(M)}{l_{\text{opt}}(M) + c_S(M)^3 + 1}$$
and that for Problem~\ref{prob:multiple} to be
$$f_2(M) = 1 + \frac{r_S(M)^2}{1+\max_S(M)+l(M)}.$$
Since both $f_1(M)$ and $f_2(M)$ depend only on the set of braids under
consideration, which doesn't change between generations, we can optimise
the simulations by caching the fitness values for a given string $M$,
rather than recalculating it each time.

\subsection{Selection}

Using an approach similar to the Stochastic Universal Sampling
Algorithm~\cite{baker:rbisa}, in each generation, we rank the candidate
sequences in order of their normalised fitness $\bar{f}(M)$: the
fitness $f(M)$ of the sequence $M$ divided by the mean fitness over the
whole population.  The integer part of $\bar{f}(M)$ gives the number
of copies contributed to the next generation, while the fractional
part gives the probability of an additional copy.  So, a sequence $M$
with a normalised fitness $\bar{f}(M) = 1.72$ contributes one copy to
the following generation, plus a 72\% chance of a second copy.

\subsection{Implementation}

The source code for the implementation is available from the authors
on request.

\section{Results}

\subsection{The single unknotting problem}

Table~\ref{tbl:braidwords} lists braid words for all knots with up to
eight crossings, and Table~\ref{tbl:results1} lists unknotting
sequences for those knots, generated by a Perl program implementing
Problem~\ref{prob:single}.  The sequences are not necessarily unique
(and in many cases will not be), nor are they guaranteed to be optimal;
however we observe that for 21 of the 35 knots with eight or fewer
crossings, our program has correctly calculated the unknotting number.

\begin{figure}
\begin{multline*}
\begin{tikzpicture}[scale=0.4]
\draw (0,0) -- (0,1);
\draw (1,0) .. controls (1,0.4) and (2,0.6) .. (2,1);
\draw[white, line width=3pt] (2,0) .. controls (2,0.4) and (1,0.6) .. (1,1);
\draw (2,0) .. controls (2,0.4) and (1,0.6) .. (1,1);
\draw (1,1) .. controls (1,1.4) and (0,1.6) .. (0,2);
\draw[white, line width=3pt] (0,1) .. controls (0,1.4) and (1,1.6) .. (1,2);
\draw (0,1) .. controls (0,1.4) and (1,1.6) .. (1,2);
\draw (2,1) -- (2,2);
\draw (0,2) -- (0,3);
\draw (1,2) .. controls (1,2.4) and (2,2.6) .. (2,3);
\draw[white, line width=3pt] (2,2) .. controls (2,2.4) and (1,2.6) .. (1,3);
\draw (2,2) .. controls (2,2.4) and (1,2.6) .. (1,3);
\draw (1,3) .. controls (1,3.4) and (0,3.6) .. (0,4);
\draw[white, line width=3pt] (0,3) .. controls (0,3.4) and (1,3.6) .. (1,4);
\draw (0,3) .. controls (0,3.4) and (1,3.6) .. (1,4);
\draw (2,3) -- (2,4);
\end{tikzpicture}
\raisebox{20pt}{$\stackrel{U^+}{\longrightarrow}$}
\begin{tikzpicture}[scale=0.4]
\draw (0,0) -- (0,1);
\draw (1,0) .. controls (1,0.4) and (2,0.6) .. (2,1);
\draw[white, line width=3pt] (2,0) .. controls (2,0.4) and (1,0.6) .. (1,1);
\draw (2,0) .. controls (2,0.4) and (1,0.6) .. (1,1);
\draw (1,1) .. controls (1,1.4) and (0,1.6) .. (0,2);
\draw[white, line width=3pt] (0,1) .. controls (0,1.4) and (1,1.6) .. (1,2);
\draw (0,1) .. controls (0,1.4) and (1,1.6) .. (1,2);
\draw (2,1) -- (2,2);
\draw (0,2) -- (0,3);
\draw (1,2) .. controls (1,2.4) and (2,2.6) .. (2,3);
\draw[white, line width=3pt] (2,2) .. controls (2,2.4) and (1,2.6) .. (1,3);
\draw (2,2) .. controls (2,2.4) and (1,2.6) .. (1,3);
\draw (0,3) .. controls (0,3.4) and (1,3.6) .. (1,4);
\draw[white, line width=3pt] (1,3) .. controls (1,3.4) and (0,3.6) .. (0,4);
\draw (1,3) .. controls (1,3.4) and (0,3.6) .. (0,4);
\draw (2,3) -- (2,4);
\end{tikzpicture}
\raisebox{20pt}{$\stackrel{M_1^-}{\longrightarrow}$}
\begin{tikzpicture}[scale=0.4]
\draw (0,0) .. controls (0,0.4) and (1,0.6) .. (1,1);
\draw[white, line width=3pt] (1,0) .. controls (1,0.4) and (0,0.6) .. (0,1);
\draw (1,0) .. controls (1,0.4) and (0,0.6) .. (0,1);
\draw (2,0) -- (2,1);
\draw (0,1) -- (0,2);
\draw (1,1) .. controls (1,1.4) and (2,1.6) .. (2,2);
\draw[white, line width=3pt] (2,1) .. controls (2,1.4) and (1,1.6) .. (1,2);
\draw (2,1) .. controls (2,1.4) and (1,1.6) .. (1,2);
\draw (1,2) .. controls (1,2.4) and (0,2.6) .. (0,3);
\draw[white, line width=3pt] (0,2) .. controls (0,2.4) and (1,2.6) .. (1,3);
\draw (0,2) .. controls (0,2.4) and (1,2.6) .. (1,3);
\draw (2,2) -- (2,3);
\draw (0,3) -- (0,4);
\draw (1,3) .. controls (1,3.4) and (2,3.6) .. (2,4);
\draw[white, line width=3pt] (2,3) .. controls (2,3.4) and (1,3.6) .. (1,4);
\draw (2,3) .. controls (2,3.4) and (1,3.6) .. (1,4);
\end{tikzpicture}
\raisebox{20pt}{$\stackrel{M_1^-}{\longrightarrow}$}
\begin{tikzpicture}[scale=0.4]
\draw (0,0) -- (0,1);
\draw (1,0) .. controls (1,0.4) and (2,0.6) .. (2,1);
\draw[white, line width=3pt] (2,0) .. controls (2,0.4) and (1,0.6) .. (1,1);
\draw (2,0) .. controls (2,0.4) and (1,0.6) .. (1,1);
\draw (0,1) .. controls (0,1.4) and (1,1.6) .. (1,2);
\draw[white, line width=3pt] (1,1) .. controls (1,1.4) and (0,1.6) .. (0,2);
\draw (1,1) .. controls (1,1.4) and (0,1.6) .. (0,2);
\draw (2,1) -- (2,2);
\draw (0,2) -- (0,3);
\draw (1,2) .. controls (1,2.4) and (2,2.6) .. (2,3);
\draw[white, line width=3pt] (2,2) .. controls (2,2.4) and (1,2.6) .. (1,3);
\draw (2,2) .. controls (2,2.4) and (1,2.6) .. (1,3);
\draw (1,3) .. controls (1,3.4) and (0,3.6) .. (0,4);
\draw[white, line width=3pt] (0,3) .. controls (0,3.4) and (1,3.6) .. (1,4);
\draw (0,3) .. controls (0,3.4) and (1,3.6) .. (1,4);
\draw (2,3) -- (2,4);
\end{tikzpicture} \\
\raisebox{20pt}{$\stackrel{\overR_3^-}{\longrightarrow}$}
\begin{tikzpicture}[scale=0.4]
\draw (0,0) .. controls (0,0.4) and (1,0.6) .. (1,1);
\draw[white, line width=3pt] (1,0) .. controls (1,0.4) and (0,0.6) .. (0,1);
\draw (1,0) .. controls (1,0.4) and (0,0.6) .. (0,1);
\draw (2,0) -- (2,1);
\draw (0,1) -- (0,2);
\draw (1,1) .. controls (1,1.4) and (2,1.6) .. (2,2);
\draw[white, line width=3pt] (2,1) .. controls (2,1.4) and (1,1.6) .. (1,2);
\draw (2,1) .. controls (2,1.4) and (1,1.6) .. (1,2);
\draw (0,2) .. controls (0,2.4) and (1,2.6) .. (1,3);
\draw[white, line width=3pt] (1,2) .. controls (1,2.4) and (0,2.6) .. (0,3);
\draw (1,2) .. controls (1,2.4) and (0,2.6) .. (0,3);
\draw (2,2) -- (2,3);
\draw (1,3) .. controls (1,3.4) and (0,3.6) .. (0,4);
\draw[white, line width=3pt] (0,3) .. controls (0,3.4) and (1,3.6) .. (1,4);
\draw (0,3) .. controls (0,3.4) and (1,3.6) .. (1,4);
\draw (2,3) -- (2,4);
\end{tikzpicture}
\raisebox{20pt}{$\stackrel{R_2^-}{\longrightarrow}$}
\begin{tikzpicture}[scale=0.4]
\draw (0,0) .. controls (0,0.4) and (1,0.6) .. (1,1);
\draw[white, line width=3pt] (1,0) .. controls (1,0.4) and (0,0.6) .. (0,1);
\draw (1,0) .. controls (1,0.4) and (0,0.6) .. (0,1);
\draw (2,0) -- (2,1);
\draw (0,1) -- (0,2);
\draw (1,1) .. controls (1,1.4) and (2,1.6) .. (2,2);
\draw[white, line width=3pt] (2,1) .. controls (2,1.4) and (1,1.6) .. (1,2);
\draw (2,1) .. controls (2,1.4) and (1,1.6) .. (1,2);
\draw (0,2) -- (0,4);
\draw (1,2) -- (1,4);
\draw (2,2) -- (2,4);
\end{tikzpicture}
\raisebox{20pt}{$\stackrel{M_2^-}{\longrightarrow}$}
\begin{tikzpicture}[scale=0.4]
\draw (0,0) .. controls (0,0.4) and (1,0.6) .. (1,1);
\draw[white, line width=3pt] (1,0) .. controls (1,0.4) and (0,0.6) .. (0,1);
\draw (1,0) .. controls (1,0.4) and (0,0.6) .. (0,1);
\draw (0,1) -- (0,4);
\draw (1,1) -- (1,4);
\draw[white] (2,0) -- (2,4);
\end{tikzpicture}
\raisebox{20pt}{$\stackrel{M_2^-}{\longrightarrow}$}
\begin{tikzpicture}[scale=0.4]
\draw (0,0) -- (0,4);
\draw[white] (1,0) -- (1,4);
\draw[white] (2,0) -- (2,4);
\end{tikzpicture}
\end{multline*}
\caption{Reduction of the figure-eight knot}
\label{fig:4_1}
\end{figure}
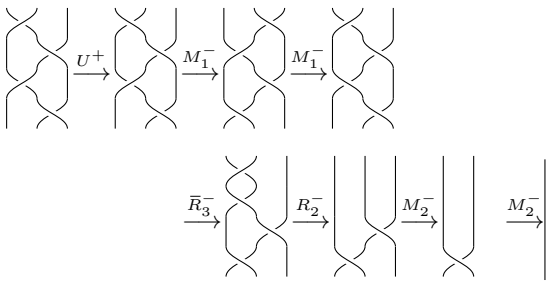

Figure~\ref{fig:4_1} shows the braid $\sigma_{1} \sigma_{2}^{-1}
\sigma_{1} \sigma_{2}^{-1}$ (whose closure is isotopic to the figure-eight
knot $4_1$) being reduced by the sequence $U M_1^2 R_3 R_2 M_2^2$.

The knots $3_1$, $5_1$ and $7_1$ are worth examining a little closer:
these are the \defn{torus knots} of type $(2,2n{+}1)$, with braid
presentation $\sigma_1^{2n+1} \in B_2$ and unknotting number $u(K)
= n$.  The unknotting sequences obtained for these knots have a very
similar form, namely $(U R_2)^{n-1} M_2$.  More generally, given two
positive, coprime integers $p$ and $q$, the \defn{torus knot} of type
$(p,q)$ is the knot which can be drawn on the surface of a standard,
unknotted torus so that the strands wind $p$ times round the torus in
the longitudinal direction, and $q$ times in the meridional direction.
The torus knot $T_{p,q}$ of type $(p,q)$ has unknotting number $u(T_{p,q})
= \tfrac12(p{-}1)(q{-}1)$.  See Rolfsen~\cite[Section~3.C]{rolfsen:kl}
or Cromwell~\cite[Section~1.5]{cromwell:ikl} for further details on
torus knots.

\begin{table}
  \centering
{\footnotesize
\begin{tabular}{|l|l|c|c|}\hline
  $K$ & Unknotting sequence $M$ & $\!c(M)\!$ & $\!u(K)\!$ \\\hline
$3_1$ & $U R_2 M_2$ & 1 & 1 \\[1ex]
$4_1$ & $U M_1^2 R_3 R_2 M_2^2$ & 1 & 1 \\[1ex]
$5_1$ & $(U R_2)^2 M_2$ & 2 & 2 \\
$5_2$ & $U R_2 M_1^2 R_3 R_2 M_2^2$ & 1 & 1 \\[1ex]
$6_1$ & $U R_2 U M_2 M_1^2 R_3 R_2 M_2^2$ & 2 & 1\\
$6_2$ & $U R_2 U M_1^2 R_3 R_2 M_2^2$ & 2 & 1\\
$6_3$ & $U R_2 U M_2 R_2 M_2$ & 2 & 1\\[1ex]
$7_1$ & $(U R_2)^3 M_2$ & 3 & 3 \\
$7_2$ & $U R_2 M_1 (M_1 U)^2 (M_1 R_3 R_2 M_2)^2 M_2$ & 3 & 1 \\
$7_3$ & $(U R_2 M_1)^2 R_3 R_2 M_2^2$ & 2 & 2 \\
$7_4$ & $U R_2 M_2 U R_2 M_1^2 R_3 R_2 M_2^2$ & 2 & 2 \\
$7_5$ & $(U R_2)^2 M_2 U R_2 M_2$ & 3 & 2 \\
$7_6$ & $U R_2 M_2 U M_1^2 R_3 R_2 M_2^2$ & 2 & 1 \\
$7_7$ & $M_1 U M_1^4 (R_3 R_2 M_2)^2 M_2$ & 1 & 1 \\[1ex]
$8_1$ & $(U R_2 M_2)^2 U M_1^2 R_3 R_2 M_2^2$ & 3 & 1 \\
$8_2$ & $U R_2 M_1^4 U M_1^2 R_3 M_1 R_2^2 M_2^2$ & 2 & 2 \\
$8_3$ & $U R_2 M_2 M_1^2 U M_1^2 R_3 R_2 M_1 M_2 R_3 R_2 M_2^2$ & 2 & 2 \\
$8_4$ & $M_1 U R_2 M_1 U M_1^3 (R_3 R_2 M_2)^2 M_2$ & 2 & 2 \\
$8_5$ & $U M_1 U R_2 M_1^2 R_3 R_2^2 M_2^2$ & 2 & 2 \\
$8_6$ & $(U R_2)^2 M_2 U M_1^2 R_3 R_2 M_2^2$ & 3 & 2 \\
$8_7$ & $M_1 U R_2 M_1 U M_1^2 R_3^2 (R_2 M_2)^2$ & 2 & 1 \\
$8_8$ & $U R_2 M_1^5 U R_2 M_1^2 R_3 M_1 R_2 M_2^3$ & 2 & 2 \\
$8_9$ & $M_1^3 U M_1^2 R_3^3 M_2 R_2^3 M_2$ & 1 & 1 \\
$8_{10}$ & $U R_2 U M_1^3 (R_3 R_2)^2 M_2^2$ & 2 & 2 \\
$8_{11}$ & $U R_2 M_2 M_1^2 U M_1^2 R_3 M_1 R_2^2 M_2^2$ & 2 & 1 \\
$8_{12}$ & $M_1 U M_1^5 U M_1 R_3 M_1 U M_2 R_2 M_2 R_3 R_2 M_2^2$ & 3 & 2 \\
$8_{13}$ & $U R_2 M_2 U R_2 M_1^2 R_3 R_2 M_2^2$ & 2 & 1 \\
$8_{14}$ & $U R_2 M_1 U M_1^2 U M_1^2 (R_3 R_2 M_2)^2 M_2$ & 3 & 1 \\
$8_{15}$ & $U R_2 M_2 M_1^3 U R_2 M_1 R_3 R_2 M_2^2$ & 2 & 2 \\
$8_{16}$ & $U R_2 U M_1 U (R_2 M_2)^2$ & 3 & 2 \\
$8_{17}$ & $M_1^2 U M_1^3 R_3 M_1 R_2 (R_3 R_2)^2 M_2^2$ & 1 & 1 \\
$8_{18}$ & $U M_1^2 R_3 R_2 M_1 U (R_2 M_2)^2$ & 2 & 2 \\
$8_{19}$ & $U R_2 M_1 R_3 U M_1 U (R_2 M_2)^2$ & 3 & 3 \\
$8_{20}$ & $M_1^5 U R_2 M_1^2 R_3 R_2^2 M_2^2$ & 1 & 1 \\
$8_{21}$ & $U R_2 M_1^3 (R_3 R_2)^2 M_2^2$ & 1 & 1 \\
\hline
\end{tabular}}
  \caption{Unknotting sequences $M$ for single knots $K$, using the braid words
  from Table~\ref{tbl:braidwords}, comparing the number $c(M)$ of
  crossing changes performed by the sequence $M$, and the crossing number
  $u(K)$ of the knot $K$}
  \label{tbl:results1}
\end{table}

\subsection{The multiple unknotting problem}

A simulation of Problem~\ref{prob:multiple}, again implemented in
Perl, obtains universal or near-universal unknotting sequences, some
examples of which are listed in Table~\ref{tbl:results2}.  Some of
the more complex braids (those corresponding to knots with minimal
crossing number 9 or higher) were unreducible by any scheme found by
our program, because we restricted ourselves to operations which don't
increase the length of the braid word.  As noted earlier (see the paper
by Coward~\cite{coward:ordering} for details), sometimes we need to
perform a move of type $\overR_2^{\pm}$ to introduce two additional
crossings during the reduction scheme.

\begin{table}
\footnotesize
\begin{tabular}{|l|l|r|r|r|} \hline
$S$             & $M$ & $\max_S(M)$ & $r_S(M)$ & $|S|$ \\ \hline
$3_1$--$4_1$    & $U M_1^2 R_3 R_2 M_2^2$     & 1 & 2 & 2 \\
$3_1$--$4_1$    & $U M_1^2 R_3 R_2 M_2$       & 2 & 2 & 2 \\
$3_1$--$4_1$    & $U R_2 R_3 M_2 M_1^2$         & 3 & 2 & 2 \\
\hline
$3_1$--$5_2$    & $M_1 U R_3 R_2 M_1 M_2^2 U$   & 2 & 4 & 4 \\
$3_1$--$5_2$    & $U R_2 M_1^2 R_3 M_2 R_2 M_2$ & 2 & 4 & 4 \\
$3_1$--$5_2$    & $U R_3 R_2 M_1^3 R_3 M_2^2$   & 2 & 4 & 4 \\
\hline
$3_1$--$6_3$    & $M_1 R_2 U M_1^2 R_3 R_2 M_2$ & 3 & 7 & 7 \\
$3_1$--$6_3$    & $M_1 U M_1^2 R_3 R_2 M_2$ & 3 & 7 & 7 \\
\hline
$3_1$--$7_7$    & $M_2 U R_2 M_2 M_1 U R_3$  & 4 & 14 & 14 \\
                & \quad $M_1^2 M_2$ & & & \\
$3_1$--$7_7$    & $M_1 R_3 U S M_1 R_2$ & 5 & 14 & 14 \\
                & \quad $M_2 M_1$ & & & \\
$3_1$--$7_7$    & $U R_3 R_2 S M_2 M_1^2 M_2$   & 6 & 14 & 14 \\
\hline
$3_1$--$8_{21}$ & $M_1 U R_3 M_1^2 R_3 M_2$ & 5 & 35 & 35 \\
                & \quad $S R_3 R_2 M_2^3$ & & & \\
$3_1$--$8_{21}$ & $M_1 U M_1 R_3 R_2 M_2^2$ & 7 & 35 & 35 \\
                & \quad $R_3 M_1$ & & & \\
\hline
$3_1$--$9_{49}$ & $U R_2 R_3 M_1^2 U M_2$ & 6 & 74 & 84 \\
                & \quad $R_3 M_1 M_2 M_1 R_3$ & & & \\
$3_1$--$9_{49}$ & $R_3 R_2 U M_1^2 M_2^2$ & 10 & 73 & 84 \\
\hline
\end{tabular}
\caption{Universal or near-universal unknotting sequences for multiple
knots}
\label{tbl:results2}
\end{table}

\section{Conclusions and future work}

In this paper we have seen how an unknotting algorithm can be evolved
based on a number of primitive moves. There are a number of future
directions for research in the area of applying evolutionary algorithms
and other machine learning techniques to mathematical problems in knot
theory and related areas.

Initially, there are some basic extensions to the work described
in this paper. For example, rather than focusing on unknotting we
could use similar techniques to address the related problem of knot
equivalence. Furthermore, there are similar problems in other areas of
mathematics (for example, graph theory and group theory) so this could
be extended to those.

More interestingly, there is the question of counterexample search. There
are a number of conjectures in this area where there are some measures
that could be used to ascertain how close a particular example is
to being a counterexample to the conjecture. Using these measures,
experiments could be done on exploring the space of knots to find a
counterexample; some preliminary work along these lines has been done
by Mahrwald~\cite{mahrwald:mmath}.

Another approach would be to take a data mining approach to certain
mathematical problems. For example, we could generate a large database
of knots and then apply classification techniques to distinguish
different classes of knot, or applying clustering techniques to group
knots according to some metric. An examination of the results from this
classification might give new insights into the underlying structure
of the space of knots. A related topic is using genetic programming to
evolve \emph{invariants}, that is, functions that distinguish between
different knots by processing their diagrams.

\bibliographystyle{AISB}
\bibliography{link}

\end{document}